\magnification=1202

\input amssym.def
%%%%%%%%%%%%%% FONTS %%%%%%%%%%%%%%%
 at9.98pt
\font\bb=msbm10 at9.98pt
 at9.98pt 
 at5pt
\font\cyr=wncyi10 at9.98pt
\font\eightrm=cmr8
\font\eightsc=cmcsc8
\font\labf=cmbx10 at13.1pt
 at15.74pt
\font\larm=cmr10 at13.1pt
 at15.74pt
\font\sc=cmcsc10 at9.98pt
\font\tenpbf=cmbx10 at8.32pt
\font\tenpit=cmti10 at8.32pt
\font\tenprm=cmr10 at8.32pt

%%%%%%%%%%%%%% SPECIAL COMMANDS %%%%%%%%%%%%%

\def\an{\raise0.5pt\hbox{$\kern2pt\scriptstyle\in\kern2pt$}}
\def\Ann{\hbox{\rm Ann\kern1pt}}
\def\Anns{\hbox{$\scriptstyle\rm Ann\kern0.5pt$}}

\def\arkef{\advance\chapternumber by 1\sc\roman{\the\chapternumber}}
\def\Aut{\hbox{\rm Aut\kern1pt}}
\def\bell{\hskip0pt\lower1.6pt\hbox{\bel\char'012}\kern5pt}

\def\callige{\hbox{\calligl e\kern2pt}}
\def\cheridexi{\hskip0pt\lower2pt\hbox{\cheridexia}\kern5pt}
\def\cheridexia{{\bbding\char'21}}

\def\coker{{\rm coker\kern1pt}}
\def\Colon{\colon\kern2pt}
\def\comp{\hbox{\lower5.8pt\hbox{\larm\char'027}}}

\def\corang{{\rm corang\kern1pt}}
\def\cos{\hbox{\rm cos\kern1pt}}
\def\cosh{\hbox{\rm cosh\kern1pt}}
\def\dbaraux{\hbox{\= {\kern-2pt\= {}}}}
\def\dbar#1{\raise3pt\hbox{\dbaraux}\kern-7.8pt #1}
\def\Der{\lower0.5pt\hbox{\ygoth Der}}

\def\dim{{\rm dim\kern1pt}}
\def\double{\hbox{\kern1.5pt\bb\char'156\kern-7.6pt\char'157\kern1.5pt}}

\def\enwsh#1{{\lower2.1pt\hbox{$\buildrel{\textstyle\cup}\over
{\lower.8pt\hbox{${}_{\scriptscriptstyle#1}$}}$}}}
\def\exp{{\rm exp\kern1pt}}

\def\exten{\hbox{\callig \kern-2.5pt Ext\lower2.5pt\hbox{\kern2.5pt}}}
\def\Ham{\hbox{\rm Ham\kern1pt}}
\def\im{{\rm im\kern1pt}}
\def\k{\raise0.25pt\hbox{$\ygot k$}}

\def\ker{{\rm ker\kern1pt}}

\def\Lie{\hbox{%$\script L$\kern-3pt\callig ie\kern2pt
\callig Lie\kern2pt}}
\def\mavrodexi{\hskip0pt\lower2pt\hbox{\mavrodexia}\kern5pt}
\def\mavrodexia{{\bbding\char'15}}
\def\meriki{\hbox{\cyr\char'144\kern0.3pt}}

\def\na{\raise0.5pt\hbox{$\kern2pt\scriptstyle\ni\kern2pt$}}
\def\noan{\hbox{$\an\raise0.6pt\hbox{$\kern-6.5pt\scriptstyle
          \slash\kern3pt$}$}}

\def\pounds{\rlap{\lower3.5pt\hbox{\kern2.9pt\hbox{\char'26}}}
           {\script L}}
\def\pr{\hbox{\kern3pt{\calligs p}\callig r\kern2pt}}
\def\qed{\hbox{\kern0.3cm\vrule height5pt width5pt depth-0.2pt}}
\def\QED{\hbox{\kern0.3cm\vrule height6pt width6pt depth-0.2pt}}
\def\R{\hbox{\bf\char'122}}
\def\rang{{\rm rang\kern1pt}}
\def\rank{{\rm rank\kern1pt}}
\def\S{\hbox{\bf\char'123}}
\def\san{\raise0.5pt\hbox{$\kern0.7pt\scriptscriptstyle
         \in\kern0.7pt$}}

\def\scomp{\hskip-0.05truecm\hbox{\lower5pt\hbox{$\mathchar"2017$}}
           \hskip-0.05truecm}

\def\sem{\hbox{{\script S}\kern-2.5pt\callig em\kern2pt}}

\def\sin{\hbox{\rm sin\kern1pt}}
\def\sinh{\hbox{\rm sinh\kern1pt}}

\def\styl{\hbox{\bbding\char'26}}
\def\stylo{\hskip0.3truecm\hbox{\lower1.5pt\hbox{\styl}}}
\def\times{\;{\mathchar"2202}\;}

\def\tonos{\hbox{\kern-1.3pt\lower0.7pt\hbox{$\mathchar"6013$}}}
\def\tonoskef{\hbox{$\kern-1.3pt\mathchar"6013$}}

\def\wbaraux{\hbox{\= {\kern-1.4pt\= {\kern-1.4pt\= {\kern-1.4pt\=
 {\kern-1.4pt\= {\kern-1.4pt\= {\kern-1.4pt\= {\kern-1.4pt\= {}}}}}}}}}}
\def\wbar#1{\hbox{\raise3pt\hbox{\wbaraux}\kern-30.5pt #1}}

\def\wwbaraux{\hbox{\= {\kern-1.4pt\= {\kern-1.4pt\= {\kern-1.4pt\=
{\kern-1.4pt\= {\kern-1.4pt\= {\kern-1.4pt\= {\kern-1.4pt\=
{\kern-1.4pt\= {}}}}}}}}}}}
\def\wwbar#1{\hbox{\raise3pt\hbox{\wwbaraux}\kern-34pt #1}}

%%%%%%%%%%%%%%%%%%% FOOTNOTE %%%%%%%%%%%%%%%%%%%
\catcode`\@=11
\def\eightpoint{\eightrm}
\def\footnote#1{\edef\@sf{\spacefactor\the\spacefactor}#1\@sf
     \insert\footins\bgroup\eightpoint
     \interlinepenalty100 \let\par=\endgraf
      \leftskip=0pt \rightskip=0pt
      \splittopskip=10pt plus 1pt minus 1pt \floatingpenalty=20000
      \smallskip\item{#1}\bgroup\strut\aftergroup\@foot\let\neft}
\skip\footins=12pt plus 2pt minus 4pt
\dimen\footins=30pc

%%%%%%%%%%%%%%%%%% CENTERLINE %%%%%%%%%%%%%%%%%%
\def\line{\hbox to\hsize}

%%%%%%%%%%%%%%%%%% TITLE %%%%%%%%%%%%%%%%%%%%
\def\title#1{\line{\hss}\line{\hss#1\hss}%
\line{\hss}\hskip-0.75truecm}

%%%%%%%%%%%%%%%%%% BIBLIOGRAPHY COMMANDS %%%%%%%%%%%%%%%%%%%%%%
\def\author#1{{\tenprm #1:}}

\def\periodiko#1{{\tenpit #1\tenprm ,}}
\def\selides#1{{\tenprm #1}}
\def\titlosa#1{{\tenprm #1,}}

\def\volume#1{{\tenprm Vol. \tenpbf #1\tenprm :}}

%
%%%%%%%%%%%%%%%%%%%%%%%%%%%%%%%%%%%%%%%%%%%%%%%%%%%%%%%%%%%%%%%%%%%%%%%%
%
%   Next, I define basic spacing parameters.
%
\def\teleia{\hbox{.}}
\newif\ifPhysRev
\def\Textindent#1{\noindent\llap{#1\enspace}\ignorespaces}
\def\nonfrenchspacing{\sfcode`\.=3001 \sfcode`\!=3000 \sfcode`\?=3000
        \sfcode`\:=2000 \sfcode`\;=1500 \sfcode`\,=1251 }
\nonfrenchspacing
\newdimen\d@twidth
 {\setbox0=\hbox{s.} \global\d@twidth=\wd0 \setbox0=\hbox{s}
        \global\advance\d@twidth by -\wd0 }
\def\removehglue{\loop \unskip \ifdim\lastskip >\z@ \repeat }
\def\roll@ver#1{\removehglue \nobreak \count255 =\spacefactor \dimen@=\z@
        \ifnum\count255 =3001 \dimen@=\d@twidth \fi
        \ifnum\count255 =1251 \dimen@=\d@twidth \fi
    \iftwelv@ \kern-\dimen@ \else \kern-0.83\dimen@ \fi
   #1\spacefactor=\count255 }
\def\step@ver#1{\relax \ifmmode #1\else \ifhmode
        \roll@ver{${}#1$}\else {\setbox0=\hbox{${}#1$}}\fi\fi }
\def\attach#1{\step@ver{\strut^{\mkern 2mu #1} }}

\normalbaselineskip = 20pt plus 0.2pt minus 0.1pt
\normallineskip = 1.5pt plus 0.1pt minus 0.1pt
\normallineskiplimit = 1.5pt
\newskip\normaldisplayskip
\normaldisplayskip = 20pt plus 5pt minus 10pt
\newskip\normaldispshortskip
\normaldispshortskip = 6pt plus 5pt
\newskip\normalparskip
\normalparskip = 6pt plus 2pt minus 1pt
\newskip\skipregister
\skipregister = 5pt plus 2pt minus 1.5pt
\newif\ifsingl@    \newif\ifdoubl@
\newif\iftwelv@    \twelv@true
\def\singlespace{\singl@true\doubl@false\spaces@t}
\def\doublespace{\singl@false\doubl@true\spaces@t}
\def\normalspace{\singl@false\doubl@false\spaces@t}
\def\Tenpoint{\tenpoint\twelv@false\spaces@t}
\def\Twelvepoint{\twelvepoint\twelv@true\spaces@t}
\def\spaces@t{\relax
      \iftwelv@ \ifsingl@\subspaces@t3:4;\else\subspaces@t1:1;\fi
       \else \ifsingl@\subspaces@t3:5;\else\subspaces@t4:5;\fi \fi
      \ifdoubl@ \multiply\baselineskip by 5
         \divide\baselineskip by 4 \fi }
\def\subspaces@t#1:#2;{
      \baselineskip = \normalbaselineskip
      \multiply\baselineskip by #1 \divide\baselineskip by #2
      \lineskip = \normallineskip
      \multiply\lineskip by #1 \divide\lineskip by #2
      \lineskiplimit = \normallineskiplimit
      \multiply\lineskiplimit by #1 \divide\lineskiplimit by #2
      \parskip = \normalparskip
      \multiply\parskip by #1 \divide\parskip by #2
      \abovedisplayskip = \normaldisplayskip
      \multiply\abovedisplayskip by #1 \divide\abovedisplayskip by #2
      \belowdisplayskip = \abovedisplayskip
      \abovedisplayshortskip = \normaldispshortskip
      \multiply\abovedisplayshortskip by #1
        \divide\abovedisplayshortskip by #2
      \belowdisplayshortskip = \abovedisplayshortskip
      \advance\belowdisplayshortskip by \belowdisplayskip
      \divide\belowdisplayshortskip by 2
      \smallskipamount = \skipregister
      \multiply\smallskipamount by #1 \divide\smallskipamount by #2
      \medskipamount = \smallskipamount \multiply\medskipamount by 2
      \bigskipamount = \smallskipamount \multiply\bigskipamount by 4 }
\def\normalbaselines{ \baselineskip=\normalbaselineskip
   \lineskip=\normallineskip \lineskiplimit=\normallineskip
   \iftwelv@\else \multiply\baselineskip by 4 \divide\baselineskip by 5
     \multiply\lineskiplimit by 4 \divide\lineskiplimit by 5
     \multiply\lineskip by 4 \divide\lineskip by 5 \fi }

%%%%%%%%%%%%%%%%%%%%%%%%%%%%%%%%%%%%%%%%%%%%%%%%%%%%%%%%%%%%%%%%%%%%%%%%
%                                  				       %
%   Here come chapter, section, subsection & appendix macros.          %
%                                                                      %
%%%%%%%%%%%%%%%%%%%%%%%%%%%%%%%%%%%%%%%%%%%%%%%%%%%%%%%%%%%%%%%%%%%%%%%%

\def\abstract#1{\parshape=1 0.7cm \dimen10
                {\tenpbf Abstract. \tenprm #1}}

\newcount\appendixnumber     \appendixnumber=0
\newcount\chapternumber      \chapternumber=0
\newcount\equanumber         \equanumber=0
\newcount\mathnumber         \mathnumber=0
\newcount\appequanumber      \appequanumber=0
\newcount\appmathnumber      \appmathnumber=0

\let\variableone=\relax
\let\variabletwo=\relax
\let\chapterlabel=\relax
\let\sectionlabel=\relax
\let\mathlabel=\relax
\newtoks\chapterstyle        \chapterstyle={\Number}
\newtoks\sectionstyle        \sectionstyle={\chapterlabel\Number}
\newskip\chapterskip         \chapterskip=\bigskipamount
\newskip\sectionskip         \sectionskip=\medskipamount
\newskip\headskip            \headskip=8pt plus 3pt minus 3pt
\newdimen\chapterminspace    \chapterminspace=15pc
\newdimen\sectionminspace    \sectionminspace=10pc
\newdimen\sectionspace       \sectionspace=20pc
\newdimen\referenceminspace  \referenceminspace=25pc

\def\chapterreset{\global\advance\chapternumber by 1
   \ifnum\equanumber<0 \else\global\equanumber=0\fi
   \mathnumber=0
   \makechapterlabel}
\def\makechapterlabel{\let\sectionlabel=\relax\let\mathlabel=\relax
 \xdef\chapterlabel{\the\chapterstyle{\the\chapternumber\teleia\kern3pt}}}

\def\rightheadline{\sc\hfil\variableone\eightsc\hfil\folio}
\def\leftheadline{\eightsc\folio\hfil{\sc\variabletwo}\hfil}
\def\heads{\footline={\hfil}\headline={\ifodd\pageno
               \rightheadline\else\leftheadline\fi}}

\def\headseis{\partreset\headline={\ifodd\pageno{
                         \hfil\sc partie {\eightsc\the\partnumber}
                         -introduction\hfil\eightsc\folio}\else
                        {\eightsc\folio\hfil\sc partie
                         {\eightsc\the\partnumber}-introduction\hfil}\fi}
                        \footline={\hfil}}

\def\alphabetic#1{\count255='140 \advance\count255 by #1\char\count255}
\def\Alphabetic#1{\count255='100 \advance\count255 by #1\char\count255}
\def\Roman#1{\uppercase\expandafter{\romannumeral #1}}
\def\roman#1{\romannumeral #1}
\def\Number#1{\number #1}
\def\BLANC#1{}

\def\titlestyle#1{\par\begingroup \interlinepenalty=9999
     \leftskip=0.02\hsize plus 0.23\hsize minus 0.02\hsize
     \rightskip=\leftskip \parfillskip=0pt
     \hyphenpenalty=9000 \exhyphenpenalty=9000
     \tolerance=9999 \pretolerance=9000
     \spaceskip=0.333em \xspaceskip=0.5em
     \iftwelv@\bf\else\bf\fi
   \noindent #1\par\endgroup }

\def\spacecheck#1{\dimen@=\pagegoal\advance\dimen@ by -\pagetotal
   \ifdim\dimen@<#1 \ifdim\dimen@>0pt \vfil\break \fi\fi}
\def\TableOfContentEntry#1#2#3{\relax}

\def\chapter#1{\par\vskip0.7cm
   %\spacecheck\chapterminspace
   \chapterreset \titlestyle{\chapterlabel\ #1}
   \nobreak\vskip\headskip
   \wlog{\string\chapter\space \chapterlabel} }

\def\appendixreset{\global\advance\appendixnumber by 1
                   \appmathnumber=0\appequanumber=0}
\def\appendix#1{\par \penalty-300\vskip\chapterskip
   \spacecheck\chapterminspace
   \appendixreset \title{\bf Appendix \Alphabetic{\the\appendixnumber}}
   \nobreak\vskip-\chapterskip\penalty 30000
   \vskip-\chapterskip
   \par{\titlestyle{#1}}
   \vskip\chapterskip
   \wlog{\string\appendix\space \chapterlabel} }

%%%%%%%%%%%%%%%%%%%%%%%%%%%%%%%%%%%%%%%%%%%%%%%%%%%%%%%%%%%%%%%%%%%%%%%%
%
%   Here come macros for equation numbering, math numbering
%
\def\eqname#1{\relax \ifnum\equanumber<0
     \xdef#1{{\noexpand\rm(\number-\equanumber)}}%
       \global\advance\equanumber by -1
    \else \global\advance\equanumber by 1
      \xdef#1{{\noexpand(
                             \rm{\number\equanumber})}} \fi #1}

\def\eqn{\eqno\eqname}

\def\math#1#2{\vskip0.1cm
   \global\advance\mathnumber by 1
   \xdef\mathlabel{\the\mathnumber}
   \wlog{\string\math\space \mathlabel}
   {\bf\enspace\mathlabel\hskip0.2cm #1}
   \xdef#2{{\mathlabel}}}

\def\appeqname#1{\relax \ifnum\appequanumber<0
     \xdef#1{{\noexpand\rm(\number-\appequanumber)}}%
       \global\advance\appequanumber by -1
    \else \global\advance\appequanumber by 1
      \xdef#1{{\noexpand(\hbox{\Alphabetic{\the\appendixnumber}}\teleia
                            {\number\appequanumber})}} \fi #1}

\def\mathapp#1#2{\vskip0.1cm
   \global\advance\appmathnumber by 1
   \xdef\appmathlabel{{\Alphabetic{\the\appendixnumber}}\teleia
   \the\appmathnumber}
   \wlog{\string\mathapp\space \appmathlabel}
   {\bf\enspace\appmathlabel\hskip0.2cm #1}
   \xdef#2{{\appmathlabel}}}

%%%%%%%%%%%%%%%%%%%%%%%%%%%%%%%%%%%%%%%%%%%%%%%%%%%%%%%%%%%%%%%%%%%%%%%%

%%%%%%%%%%%%%%%%%%%%%%%%%%%%%%%%%%%%%%%%%%%%%%%%%%%%%%%%%%%%%%%%%%%%%%%%
%
%   Here come macros for references, figures & tables.
%
% % % % % % % % % % % % % % % % % % % % % % % % % % % % % % % % % % % %
%%  First, references.
%
\newtoks\referencestyle      \referencestyle={\tenpbf\Number}
\newcount\referencecount     \referencecount=0
\newcount\lastrefsbegincount \lastrefsbegincount=0
\newif\ifreferenceopen       \newwrite\referencewrite
\newif\ifrw@trailer
\newdimen\refindent     \refindent=13pt
\def\NPrefmark#1{\attach{\scriptscriptstyle [ #1 ] }}
\let\PRrefmark=\attach
\def\refmark#1{\relax\ifPhysRev\PRrefmark{#1}\else\NPrefmark{#1}\fi}
\def\refend@{\refmark{\number\referencecount}}
\def\refend{\refend@{}\space }
\def\refsend{\refmark{\count255=\referencecount
   \advance\count255 by-\lastrefsbegincount
   \ifcase\count255 \number\referencecount
   \or \number\lastrefsbegincount,\number\referencecount
   \else \number\lastrefsbegincount-\number\referencecount \fi}\space }
\def\refitem#1{\par\hangafter=0 \hangindent=\refindent	\Textindent{#1}}
\def\Ref{\rw@trailertrue\REF}
\def\REF#1{\r@fstart{#1}%
   \rw@begin{\tenprm [\tenpbf\Number{\the\referencecount}\tenprm ]}\rw@end}
\def\r@fstart#1{\chardef\rw@write=\referencewrite \let\rw@ending=\refend@
   \ifreferenceopen \else \global\referenceopentrue
   \immediate\openout\referencewrite=referenc.txa
   \toks0={\catcode`\^^M=10}\immediate\write\rw@write{\the\toks0} \fi
   \global\advance\referencecount by 1 
   \xdef#1{[{\the\referencestyle{\the\referencecount}}]}}
 {\catcode`\^^M=\active %
 \gdef\rw@begin#1{\immediate\write\rw@write{\noexpand\refitem{#1}}%
   \begingroup \catcode`\^^M=\active \let^^M=\relax}%
 \gdef\rw@end#1{\rw@@end #1^^M\rw@terminate \endgroup%
   \ifrw@trailer\rw@ending\global\rw@trailerfalse\fi }%
 \gdef\rw@@end#1^^M{\toks0={#1}\immediate\write\rw@write{\the\toks0}%
   \futurelet\n@xt\rw@test}%
 \gdef\rw@test{\ifx\n@xt\rw@terminate \let\n@xt=\relax%
       \else \let\n@xt=\rw@@end \fi \n@xt}%
}
\let\rw@ending=\relax
\let\rw@terminate=\relax

\def\Closeout#1{\toks0={\catcode`\^^M=5}\immediate\write#1{\the\toks0}%
   \immediate\closeout#1}

\topskip1truecm
\voffset=2.5truecm
\hsize 15truecm
\vsize 20truecm
\hoffset=0.5truecm
\def\undertext#1{$\underline{\hbox{#1}}$}
\topglue 3truecm
\nopagenumbers
\def\variableone{p. baguis, m. cahen}
\def\variabletwo{a construction of symplectic connections through reduction}
\def\zero{\mathaccent"0017}
\def\absize{11cm}
\def\abstract#1{\baselineskip=12pt plus .2pt
                \parshape=1 0.7cm \absize
                 {\tenpbf Abstract. \tenprm #1}}
%%%%%%%%%%%%%%%%%%%%%%%%%%%%%%%%%%%%%%%%%
%%                 REFERENCES              %%
%%%%%%%%%%%%%%%%%%%%%%%%%%%%%%%%%%%%%%%%%
\Ref\devischer{
\author{M. DeVischer}
\titlosa{M\'emoire de Licence}
\selides{Bruxelles, 2000}}
				 
\Ref\vaisman{
\author{I. Vaisman}
\titlosa{Symplectic curvature tensors}
\periodiko{Mon. Math.}
\volume{100}
\selides{299--327 (1985)}}

\Ref\cahen{
\author{M. Cahen,  S. Gutt, J. Rawnsley}
\titlosa{Symmetric symplectic spaces with Ricci type curvature}
\periodiko{Math. Phys. Studies}
\volume{22(2)}
\selides{81--93 (2000)}}

%%%%%%%%%%%%%%%%%%%%%%%%%%%%%%%%%%%%%%%%%

\hskip10cm

\centerline{\labf A construction of symplectic connections}
\vskip0.2cm
\centerline{\labf through reduction}

\vskip1cm

\centerline{\bf P. Baguis\footnote{$^{1}$}{e-mail: 
pbaguis@ulb.ac.be}$^{,}$\footnote{$^{2}$}{Research 
supported by the Marie Curie Fellowship Nr. 
HPMF-CT-1999-00062 }, M. Cahen\footnote{$^{3}$}{e-mail: 
mcahen@ulb.ac.be}$^{,}$\footnote{$^{4}$}{Research 
supported by an ARC of the ``Communaut\'e fran\c caise de Belgique''}}

\vskip0.3cm

\centerline{Universit\'e Libre de Bruxelles}
\centerline{Campus Plaine, CP 218 Bd du Triomphe}
\centerline{1050, Brussels, Belgium}

\vskip1cm

\abstract{We give an elementary construction of symplectic 
connections through reduction. This provides an elegant description 
of a class of symmetric spaces and gives examples of symplectic 
connections with Ricci type curvature, which are not locally symmetric; 
the existence of such symplectic connections was unknown.}

\vskip1cm

{\tenprm 
{\tenpit Key-words}: Marsden-Weinstein reduction, symplectic 
connections, symmetric spaces

\vskip0.2cm

{\tenpit MSC 2000}: 53B05, 53D20, 53C15}

\vfill\eject

\baselineskip=14pt plus .2pt

\heads

{\bf 1.} Let $(M,\omega)$ be a smooth $2n$-dimensional symplectic manifold; 
a linear connection $\nabla$ on $(M,\omega)$ is said to be symplectic 
if it is torsion free and if $\omega$ is parallel. If $n>1$ the 
curvature tensor $R$  of such a connection has two irreducible 
components under the pointwise action of the linear symplectic group
$Sp(n,{\bf R})$ {\devischer} {\vaisman}. We shall denote them by $E$ 
and $W$:
$$R=E+W.\eqn\one$$
The $E$ component is determined by the Ricci tensor of $\nabla$; if 
the $W$ component vanishes the curvature is said to be of Ricci type.

In {\cahen} the simply connected symmetric symplectic spaces, whose 
curvature is of Ricci type have been classified algebraically. It was 
shown that an isomorphism class was determined by the orbit, under 
the action of the linear symplectic group $Sp(n,{\bf R})$, of an 
element $A$ belonging to the Lie algebra $\frak s\frak p(n,{\bf R})$
of the linear symplectic group such that 
$$A^{2}=\lambda\;{\rm Id}.\eqn\two$$
where $\lambda$ is any real number. It was also observed that the only 
compact symmetric symplectic space with non-zero curvature of 
Ricci type is the complex projective space ${\bf P}_{n}({\bf C})$. 

In this paper we give a very elementary and geometrical construction 
of those symmetric symplectic spaces and we provide examples of 
symplectic connections with Ricci type curvature, which are not 
locally symmetric. Finally we give a suggestion how to generalize 
this easy construction.

\vskip0.5cm

{\bf 2.} Let $A\neq 0$ be an element of $\frak s\frak p(n+1,{\bf R})$
and let  $H$ be a homogeneous polynomial of degree  2 on ${\bf R}^{2n+2}$
defined by:
$$H(x)=\omega(x,Ax), \quad\forall x\an{\bf R}^{2n+2}\eqn\three$$
where $\omega$ is the standard symplectic structure on ${\bf R}^{2n+2}$.
Let $0\neq\mu_{0}\an{\bf R}$ and denote by $\Sigma_{\mu_{0}}$ the quadric 
on ${\bf R}^{2n+2}$:
$$\Sigma_{\mu_{0}}=\big\{x\an{\bf 
R}^{2n+2}\;|\;H(x)=\mu_{0}\big\}.\eqn\four$$
This is a closed embedded submanifold of ${\bf R}^{2n+2}$. The tangent 
space to $\Sigma_{\mu_{0}}$ at $x$ is $(Ax)^{\perp}$ (with respect 
to $\omega$); the restriction of $\omega$ to $\Sigma_{\mu_{0}}$ admits 
at the point $x$ a 1-dimensional radical which is spanned by the 
Hamiltonian vector field
$$X_{H}(x)=-2Ax.\eqn\five$$
Observe that the condition $\mu_{0}\neq 0$, implies that the radial 
vector $x$ is, at the point $x$, transversal to $\Sigma_{\mu_{0}}$.
 
Let $\zero\nabla$ be the standard flat affine 
connection on ${\bf R}^{2n+2}$; it is clearly symplectic relative 
to $\omega$. Let $Y,Z$ be smooth vector fields on $\Sigma_{\mu_{0}}$.
Define a linear connection $\nabla$ on $\Sigma_{\mu_{0}}$ by:
$$(\nabla_{Y}Z)_{x}=(\zero\nabla_{Y}Z)_{x}+{1\over\mu_{0}}\omega(Z,AY)x.
      \eqn\six$$
This is indeed a vector belonging to the tangent space 
$T_{x}\Sigma_{\mu_{0}}$:
$$\eqalign{\omega_{x}(\nabla_{Y}Z,Ax)&=\omega_{x}(\zero\nabla_{Y}Z,Ax)+
\omega(Z,AY)\cr
\hfill&=-\omega(Z,AY)+\omega(Z,AY)=0.\cr}$$
This connection is torsion free as:
$$\eqalign{\nabla_{Y}Z-\nabla_{Z}Y&=\zero\nabla_{Y}Z-\zero\nabla_{Z}Y+
{1\over\mu_{0}}\big((\omega(Z,AY)-\omega(Y,AZ)\big)x\cr
\hfill&=[Y,Z].\cr}$$

\math{Lemma.}{\tangentgeodesics}{\sl The Hamiltonian vector field 
$X_{H}$, restricted to $\Sigma_{\mu_{0}}$, is tangent to geodesics of
$\nabla$ if and only if $$A^{2}=\lambda \;{\rm Id}.$$}

\undertext{\it Proof.} Applying the definition:
$${1\over 4}\nabla_{X_{H}}X_{H}(x)=A^{2}x+
{1\over\mu_{0}}\omega(Ax,A^{2}x)x.$$
If $X_{H}$ is tangent to geodesics of $\nabla$, there exists a 
function $\nu(x)$ such that
$$A^{2}x+{1\over\mu_{0}}\omega(Ax,A^{2}x)x=\nu(x)Ax.$$
But this implies $\nu(x)=0$. Deriving this relation in any direction 
$Y$ tangent to $\Sigma_{\mu_{0}}$:
$$ A^{2}Y+{2\over\mu_{0}}\omega(AY,A^{2}x)x+
{1\over\mu_{0}}\omega(Ax,A^{2}x)Y=0.$$
But using the relation once more:
$$\omega(AY,A^{2}x)=0.$$ Hence
$$A^{2}=\lambda \;{\rm Id}.$$ The converse is obvious.\qed

\vskip0.3cm

\math{Lemma.}{\parametricgroup}{\sl The Hamiltonian vector field 
$X_{H}$ generates a one parametric group $\psi_{t}$ of affine 
transformations of $(\Sigma_{\mu_{0}},\nabla)$.}

\undertext{\it Proof.} Let $Y,Z$ be smooth vector fields along 
$\Sigma_{\mu_{0}}$. Then:
$$\eqalign{[X_{H},\nabla_{Y}Z]&=[X_{H},\zero\nabla_{Y}Z+
{1\over\mu_{0}}\omega(Z,AY)x]\cr
\hfill&=\zero\nabla_{[X_{H},Y]}Z+\zero\nabla_{Y}[X_{H},Z]+\cr
\hfill&\quad{1\over\mu_{0}}\big(\omega([X_{H},Z],AY)+\omega(Z,[X_{H},AY])\big)x\cr}$$
as $X_{H}$ is an affine vector field for $\zero\nabla$.

On the other hand:
$$\eqalign{\nabla_{[X_{H},Y]}Z+\nabla_{Y}[X_{H},Z]&=\zero\nabla_{[X_{H},Y]}Z+
{1\over\mu_{0}}\omega(Z,A[X_{H},Y])x\cr
\hfill&\kern12pt+\zero\nabla_{Y}[X_{H},Z]+{1\over\mu_{0}}\omega([X_{H},Z],AY)x.\cr}$$
Hence the conclusion.\qed

\vskip0.5cm

{\bf 3.} The action of $\psi_{t}$ on $\Sigma_{\mu_{0}}$ is free. We 
describe, when $A^{2}=\lambda \;{\rm Id}$, $A\neq 0$,
the orbit space $M=\Sigma_{\mu_{0}}/\psi_{t}$.
For all values of $\lambda$, $M$ is a smooth manifold and the canonical projection 
$\pi:\Sigma_{\mu_{0}}\rightarrow M$ is a smooth submersion.

\vskip0.3cm

{\sl Case 1.} Assume  $A^{2}=-{\rm Id}$. Then:
$\omega(Ax,Ay)=-\omega(x,A^{2}y)=\omega(x,y)$
and thus $A\an Sp(n+1,{\bf R})$. There exists $0\neq x\an{\bf 
R}^{2n+2}$ such that $$\omega(x,Ax)=\epsilon_{1}, 
\quad\epsilon_{1}^{2}=1.\eqn\seven$$
Then:
$$\omega(Ax,A^{2}x)=\epsilon_{1}.\eqn\eignt$$
Hence, by recurrence, there exists a basis $\{e_{a};a\leq 2n+2\}$ 
of ${\bf R}^{2n+2}$ and an integer $p$ ($0\leq p\leq n+1$) such that
$$\omega(e_{2i-1},e_{2i}=Ae_{2i-1})=1, \quad i\leq p\eqn\nine$$
$$\omega(e_{2j-1},e_{2j}=Ae_{2j-1})=-1, \quad  p<j\leq n+1.\eqn\ten$$
We shall choose $\mu_{0}=\pm 1$ and denote $q=n+1-p$. We have the 
obvious isomorphism
$$\Sigma_{p,q,1}=\Sigma_{q,p,-1},\eqn\eleven$$
where $\Sigma_{p,q,\epsilon}$ denotes the quadric
$$\sum_{i=1}^{p}[(x^{2i-1})^{2}+(x^{2i})^{2}]-
\sum_{j=p+1}^{n+1}[(x^{2j-1})^{2}+(x^{2j})^{2}]=\epsilon.\eqn\twelve$$
The following are obvious:
$$\Sigma_{n+1,0,1}={\bf S}^{2n+1}\eqn\thirteen$$
$$\Sigma_{p,q,1}={\bf S}^{2p-1}\times{\bf R}^{2q}\quad n+1>p\geq 
1.\eqn\fourteen$$
The Hamiltonian vector field $X_{H}$ generates an action of $U(1)$ and 
one easily checks the

\math{Lemma.}{\reducedmanifold}{\sl The reduced manifold 
$\Sigma_{p,q,1}/U(1)$ is one of the following ones:
$$\Sigma_{n+1,0,1}/U(1)={\bf P}_{n}({\bf C})$$
$$\Sigma_{1,n,1}/U(1)={\bf C}^{n}.$$
For $n>p\geq 2$ $$\Sigma_{p,q,1}/U(1)$$ is a rank $q$ complex vector 
bundle over ${\bf P}_{p-1}({\bf C})$.}

\vskip0.3cm

{\sl Case 2.} Assume $A^{2}=1$. Let $V^{\pm}=\big\{x\an{\bf 
R}^{2n+2}\;|\;Ax=\pm x\big\}$. Then, if $x,y\an V^{\pm}$, 
$\omega(x,y)=\omega(Ax,Ay)=-\omega(x,A^{2}y)=-\omega(x,y)=0$.
Hence $V^{+}$ and $V^{-}$ are two supplementary lagrangian subspaces.
Choose a basis $\{e_{i};i\leq n+1\}$ of $V^{+}$ and a basis 
$\{f_{j};j\leq n+1\}$ of $V^{-}$ such that 
$$\omega(e_{i},f_{j})=\delta_{ij}.\eqn\fifteen$$
Then, if $\mu_{0}=-2$,
$$\Sigma_{\mu_{0}}=\Bigg\{x=\sum_{i=1}^{n+1}x^{i}e_{i}+y^{i}f_{i}\;\Big|\;
\sum_{i=1}^{n+1}x^{i}y^{i}=1\Bigg\}.\eqn\sixteen$$
The Hamiltonian vector field
$$-{1\over 2} X_{H}=\sum_{i=1}^{n+1}x^{i}e_{i}-y^{i}f_{i}\eqn\seventeen$$
has trajectories:
$$x^{i}(t)=x^{i}(0)e^{t}\eqn\eighteen$$
$$y^{i}(t)=y^{i}(0)e^{-t}.\eqn\nineteen$$
On the orbit of $(x(0),y(0))$ there is a unique point $(\bar{x},\bar{y})$ 
such that $\sum_{i=1}^{n+1}(\bar{x}^{i})^{2}=1$. Define then the 
point $z\an{\bf R}^{n+1}$ by
$$\bar{y}^{i}=\bar{x}^{i}+z^{i}.\eqn\twenty$$
Then:
$$\sum_{i=1}^{n+1}\bar{x}^{i}z^{i}=0.\eqn\twentyone$$
Hence we have:

\math{Lemma.}{\isomorphic}{\sl The reduced manifold $\Sigma_{-2}/{\bf R}$
is isomorphic to the tangent bundle $T{\bf S}^{n}$.}

\vskip0.3cm

{\sl Case 3.} Assume $A^{2}=0$ ($A\neq 0$). Let $V=\im A$ and $W=\ker A$. 
Then $V\subset W$.  Furthermore $(\im 
A)^{\perp}=\big\{y\;|\;\omega(Ax,y)=0,\forall x\big\}$. Hence $\omega(x,Ay)=0$, 
$\forall x$; hence $Ay=0$ and $(\im A)^{\perp}\subset\ker A$. Hence 
by dimension 
$$(\im A)^{\perp}=\ker A.\eqn\twentytwo$$
Let $X$ be a supplementary subspace to $V$ in $W$; then $X$ is 
symplectic. Let $1\leq p\leq n+1$ be the dimension of $V$. Let
$\{u_{1},\ldots,u_{n+1-p},v_{1},\ldots,v_{n+1-p}\}$ be a basis of $X$
such that:
$$\omega(u_{i},v_{j})=\delta_{ij}\eqn\twentythree$$
$$\omega(u_{i},u_{j})=\omega(v_{i},v_{j})=0\eqn\twentyfour$$
The subspace $X^{\perp}$ is symplectic and contains $V$ as a 
lagrangian subspace. Let $\{e_{k};k\leq p\}$ be a basis of a 
lagrangian subspace $V^{\star}$ of $X^{\perp}$, supplementary to $V$, 
such that:
$$\omega(e_{k},Ae_{\ell})=\epsilon_{k}\delta_{k\ell},\eqn\twentyfive$$
where $\epsilon_{k}=1$ if $k\leq q$, and $\epsilon_{k}=-1$ if $q<k$. Choose 
$\mu_{0}=\epsilon$ $(\epsilon^{2}=1)$. Then 
$$\Sigma_{p,q,\epsilon}=\Bigg\{x\an{\bf R}^{2n+2}\;\Big|\;
\sum_{k=1}^{p}\epsilon_{k}(x^{k})^{2}=\epsilon\Bigg\}.\eqn\twentysix$$
We have the isomorphism $\Sigma_{p,q,\epsilon}=\Sigma_{q,p,-\epsilon}$.
If $q=0$, $\epsilon$ must be chosen to be -1. Clearly,
$$\Sigma_{1,1,1}=2 \,{\rm points}\times{\bf R}^{2n+1}\eqn\twentyseven$$
$$\Sigma_{p,p,1}={\bf S}^{p-1}\times{\bf R}^{2n+2-p},\quad 
p>1\eqn\twentyeight$$
$$\Sigma_{p,1,1}=({\bf R}^{p-1}\cup{\bf R}^{p-1})\times{\bf 
R}^{2n+2-p},
\quad p>1\eqn\twentynine$$
$$\Sigma_{p,q,1}=({\bf S}^{q-1}\times{\bf R}^{p-q})\times{\bf 
R}^{2n+2-p},
\quad p>1,q>1.\eqn\thirty$$
The Hamiltonian vector field reads
$$X_{H}=-2\sum_{k=1}^{p}x^{k}Ae_{k}.\eqn\thirtyone$$
It generates an action of ${\bf R}$; if
$$x=\sum_{k=1}^{p}x^{k}e_{k}+y^{k}Ae_{k}+
          \sum_{a=1}^{n+1-p}z^{a}u_{a}+(z^{\prime})^{a}v_{a},\eqn\thirtytwo$$
we have:
$$x^{k}(t)=x^{k}(0)\quad z^{a}(t)=z^{a}(0)\quad 
(z^{\prime})^{a}(t)=(z^{\prime})^{a}(0)\eqn\thirtythree$$
$$y^{k}(t)=y^{k}(0)-2x^{k}(0)t.\eqn\thirtyfour$$
Observe that, if $\mu_{0}=1$,
$$\sum_{k=1}^{p}\epsilon_{k}x^{k}(t)y^{k}(t)=
-2t+\sum_{k=1}^{p}\epsilon_{k}x^{k}(0)y^{k}(0).\eqn\thirtyfive$$
Hence there exists a unique point on the orbit such that 
$\sum_{k=1}^{p}\epsilon_{k}x^{k}y^{k}=0$. Thus:

\math{Lemma.}{\reduced}{\sl The reduced manifold $\Sigma_{p,q,1}/{\bf R}$
is one of the following ones:

$$\Sigma_{1,1,1}/{\bf R}={\bf R}^{2n}\cup{\bf R}^{2n}$$
$$\Sigma_{p,1,1}/{\bf R}=({\bf R}^{p-1}\cup{\bf R}^{p-1})
\times{\bf R}^{p-1}\times{\bf R}^{2n+2-2p}\quad p>1$$
$$\Sigma_{p,q,1}/{\bf R}=T({\bf S}^{q-1}\times{\bf R}^{p-q})
\times{\bf R}^{2n+2-2p}\quad p>1,q>1$$}

\vskip0.3cm

{\bf Remark.} The condition $A^{2}=\lambda \;{\rm Id}$ is not a 
necessary condition for the quotient $\Sigma_{\mu_{0}}/\psi_{t}$
to have a natural structure of manifold such that the projection be a 
smooth submersion. Consider $A\an\frak s\frak p(2,\R)$ having a
complex eigenvalue $\lambda=a+ib$ ($ab\neq 0$). Then $A$ admits as 
other eigenvalues $\bar\lambda,-\lambda,-\bar\lambda$. Let $\{e_{\mu}\}$ 
be a basis of the complexified eigenspace corresponding to $\mu$. 
Then, we can choose the $e_{\mu}$'s such that:
$$Ae_{\mu}=\mu e_{\mu},$$
$$\omega(e_{\lambda},e_{\bar \lambda})=\omega(e_{\lambda},e_{-\bar \lambda})=
\omega(e_{\bar\lambda},e_{-\lambda})=\omega(e_{-\lambda},e_{-\bar\lambda})=0,$$
$$\overline{\omega(e_{\lambda},e_{-\lambda})}=\omega(e_{\bar\lambda},
e_{-\bar\lambda})\neq 0.$$
This complex basis is determined up to a transformation
$$e_{\lambda}\mapsto\rho e_{\lambda},\quad e_{-\lambda}\mapsto\sigma 
e_{-\lambda},\quad \rho\sigma\neq 0.$$
Choose the factors in such a way that:
$$\omega(e_{\lambda},e_{-\lambda})=\omega(e_{\bar\lambda},
e_{-\bar\lambda})=1.$$
Then writing
$$e_{\lambda}=e_{1}+ie_{2}\quad e_{-\lambda}=e_{3}+ie_{4}$$
we see that:
$$Ae_{1}=ae_{1}-be_{2}$$
$$Ae_{2}=be_{1}+ae_{2}$$
$$Ae_{3}=-ae_{3}+be_{4}$$
$$Ae_{4}=-be_{3}-ae_{4}$$
$$\omega(e_{1},e_{2})=\omega(e_{2},e_{3})=\omega(e_{1},e_{4})
=\omega(e_{3},e_{4})=0$$
$$\omega(e_{1},e_{3})={1\over 2}=-\omega(e_{2},e_{4}).$$

The quadric $H(x)=\omega(x,Ax)=1$ reads:
$$-x^{1}(ax^{3}+bx^{4})+x^{2}(ax^{4}-bx^{3})=1$$
and is thus diffeomorphic to $\S^{1}\times\R^{2}$. The orbits of the 
Hamiltonian vector field are
$$x^{1}+ix^{2}=(x^{1}+ix^{2})(0)\exp((a-ib)t)$$
$$x^{3}+ix^{4}=(x^{3}+ix^{4})(0)\exp(-(a-ib)t).$$
We can rewrite the equation $H(x)=1$ as:
$$\frak R\frak e[(x^{1}+ix^{2})(x^{3}+ix^{4})(-a+ib)]=1.$$
Hence on each orbit there exists a unique point such that 
$|x^{1}+ix^{2}|=1$. Hence the quotient is a cylinder $\S^{1}\times\R$.

\vskip0.5cm

{\bf 4.} Let $\pi: \Sigma_{\mu_{0}}\rightarrow M=\Sigma_{\mu_{0}}/\psi_{t}$
be the canonical projection and let $y=\pi(x)$. Let $H_{x}$ be the 
subspace of the tangent space to $\Sigma_{\mu_{0}}$ at $x$ which is
$H_{x}=[{\rm span}(x,Ax)]^{\perp}$. The differential of $\pi$ at $x$, 
$\pi_{\ast x}$, is a 
linear isomorphism $H_{x}\rightarrow M_{y}$. Hence a vector field $Y$ 
on $M$ admits a unique lift $\bar{Y}$ to $\Sigma_{\mu_{0}}$, which 
belongs at each point $x$, to $H_{x}$.

Notice that
$$\omega_{x}(\nabla_{\bar{Y}}\bar{Z},x)=
\omega_{x}(\zero\nabla_{\bar{Y}}\bar{Z},x)=
-\omega_{x}(\bar{Z},\bar{Y}).\eqn\thirtysix$$
Hence define the reduced connection $\nabla^{r}$ on $M$ by:
$$\overline{\nabla^{r}_{Y}Z}=\nabla_{\bar{Y}}\bar{Z}+
{1\over{\mu_{0}}}\omega(\bar{Y},\bar{Z})Ax.\eqn\thirtyseven$$

\math{Lemma.}{\reduced}{\sl The reduced connection $\nabla^{r}$ 
on $M$ is symplectic with respect to the symplectic form $\Omega$:
$$\Omega_{y}(Y,Z)=\omega_{x}(\bar{Y},\bar{Z}).\eqn\thirtyeight$$}
\undertext{\it Proof.} The torsion free condition reads:
$$\eqalign{0&=\overline{\nabla^{r}_{Y}Z}-\overline{\nabla^{r}_{Z}Y}-\overline{[Y,Z]}\cr
\hfill&=\nabla_{\bar{Y}}\bar{Z}+{1\over{\mu_{0}}}\omega(\bar{Y},\bar{Z})Ax-
\nabla_{\bar{Z}}\bar{Y}+{1\over{\mu_{0}}}\omega(\bar{Z},\bar{Y})Ax-\overline{[Y,Z]}.\cr}$$
Now:
$$\pi_{\ast}\overline{[Y,Z]}=\pi_{\ast}{[\bar{Y},\bar{Z}]}=[Y,Z]$$ and
$$\eqalign{\omega([\bar{Y},\bar{Z}],x)&= 
                        \omega(\nabla_{\bar{Y}}\bar{Z}-\nabla_{\bar{Z}}\bar{Y},x)\cr
                        \hfill&=-\omega(\bar{Z},\bar{Y})+\omega(\bar{Y},\bar{Z})\cr}$$
i.e.
$$\overline{[Y,Z]}=[\bar{Y},\bar{Z}]+{2\over\mu_{0}}\omega(\bar{Y},\bar{Z})Ax.$$
Hence the torsion free condition is satisfied. To prove that $\Omega$ 
is parallel we note that:
$$\eqalign{Y\Omega(Z,U)&= \bar{Y}\omega({\bar{Z}},\bar{U})=
\omega(\zero\nabla_{\bar{Y}}\bar{Z},\bar{U})+
\omega(\bar{Z},\zero\nabla_{\bar{Y}}\bar{U})\cr
\hfill&=\omega(\overline{\nabla^{r}_{Y}Z},\bar{U})+
                  \omega(\bar{Z},\overline{\nabla^{r}_{Y}U})\cr
\hfill&=\Omega(\nabla^{r}_{Y}Z,U)+
                    \Omega(Z,\nabla^{r}_{Y}U).\cr}$$\qed
\vskip0.3cm

\math{Theorem.}{\redcurv}{\sl The curvature tensor $R^{r}$ of the 
reduced connection $\nabla^{r}$ on $(M,\Omega)$ is of Ricci type.}

\undertext{\it Proof.} The lift of the
curvature endomorphism of $\nabla^{r}$ is 
given by:
\nobreak
$$\overline{R^{r}(Y,Z)T}=\overline{\nabla^{r}_{Y}\nabla^{r}_{Z}T}-
\overline{\nabla^{r}_{Z}\nabla^{r}_{Y}T}-\overline{\nabla^{r}_{[Y,Z]}T}.$$
\nobreak
We have:
$$\eqalign{\overline{\nabla^{r}_{Y}\nabla^{r}_{Z}T}&=\nabla_{\bar 
Y}\overline{\nabla^{r}_{Z}T}+{1\over\mu_{0}}\omega(\bar 
Y,\overline{\nabla^{r}_{Z}T})Ax\cr
\hfill&=\nabla_{\bar Y}(\nabla_{\bar Z}\bar 
T+{1\over\mu_{0}}\omega(\bar Z,\bar T)Ax)+{1\over\mu_{0}}\omega(\bar 
Y,\nabla_{\bar Z}\bar T)Ax\cr
\hfill&=\nabla_{\bar Y}\nabla_{\bar Z}\bar T+{1\over\mu_{0}}\bar Y\omega
(\bar Z,\bar T)Ax+{1\over\mu_{0}}\omega(\bar Z,\bar T)\nabla_{\bar 
Y}Ax+{1\over\mu_{0}}\omega(\bar 
Y,\nabla_{\bar Z}\bar T)Ax\cr
\hfill&=\zero\nabla_{\bar Y}(\zero\nabla_{\bar Z}\bar 
T+{1\over\mu_{0}}\omega(\bar T,A\bar 
Z)x)+{1\over\mu_{0}}\omega(\zero\nabla_{\bar Z}\bar 
T+{1\over\mu_{0}}\omega(\bar T,A\bar Z)x,A\bar 
Y)x\cr
\hfill&\kern12pt+{1\over\mu_{0}}[\omega(\zero\nabla_{\bar Y}\bar Z,\bar 
T)+\omega(\bar Z,\zero\nabla_{\bar Y}\bar 
T)]Ax\cr
\hfill&\kern12pt+{1\over\mu_{0}}\omega(\bar Z,\bar T)[A\bar 
Y+{1\over\mu_{0}}\omega(Ax,A\bar Y)x]\cr
\hfill&\kern12pt+{1\over\mu_{0}}\omega(\bar 
Y,\zero\nabla_{\bar Z}\bar T+{1\over\mu_{0}}\omega(\bar T,A\bar 
Z)x)Ax\cr
&=\zero\nabla_{\bar Y}\zero\nabla_{\bar Z}\bar 
T+{1\over\mu_{0}}[\omega(\zero\nabla_{\bar Y}\bar T,A\bar 
Z)+\omega(\bar T,\zero\nabla_{\bar Y}A\bar Z)]x\cr
\hfill&\kern12pt+{1\over\mu_{0}}\omega(\bar T,A\bar Z)\bar 
Y+{1\over\mu_{0}}\omega(\zero\nabla_{\bar Z}\bar T,A\bar Y)x\cr
\hfill&\kern12pt+{1\over\mu_{0}}[\omega(\zero\nabla_{\bar Y}\bar 
Z,\bar T)+\omega(\bar Z,\zero\nabla_{\bar Y}\bar T)]Ax\cr
\hfill&\kern12pt+{1\over\mu_{0}}\omega(\bar Z,\bar T)A\bar 
Y+{1\over\mu_{0}^{2}}\omega(\bar Z,\bar T)\omega(Ax,A\bar 
Y)x+{1\over\mu_{0}}\omega(\bar Y,\zero\nabla_{\bar Z}\bar T)Ax\cr}$$
$$\eqalign{
\hfill&=\zero\nabla_{\bar Y}\zero\nabla_{\bar Z}\bar 
T+{1\over\mu_{0}}\omega(\bar T,A\bar Z)\bar 
Y+{1\over\mu_{0}}\omega(\bar Z,\bar T)A\bar Y\cr
\hfill&\kern12pt+{1\over\mu_{0}}[\omega(\zero\nabla_{\bar Y}\bar 
T,A\bar Z)+\omega(\bar T,\zero\nabla_{\bar Y}A\bar 
Z)\cr
\hfill&\kern38pt+\omega(\zero\nabla_{\bar Z}\bar T,A\bar 
Y)+{1\over\mu_{0}}\omega(\bar Z,\bar T)\omega(Ax,A\bar Y)]x\cr
\hfill&\kern12pt+{1\over\mu_{0}}[\omega(\zero\nabla_{\bar Y}\bar 
Z,\bar T)+\omega(\bar Z,\zero\nabla_{\bar Y}\bar T)+\omega(\bar 
Y,\zero\nabla_{\bar Z}\bar T)]Ax.}$$
Also:
$$\eqalign{\overline{\nabla^{r}_{[Y,Z]}T}&=\nabla_{[Y,Z]}\bar 
T+{1\over\mu_{0}}\omega(\overline{[Y,Z]},\bar T)Ax\cr
\hfill&=\zero\nabla_{\overline{[Y,Z]}}\bar 
T+{1\over\mu_{0}}\omega(\bar 
T,A\overline{[Y,Z]})x+{1\over\mu_{0}}\omega(\overline{[Y,Z]},\bar 
T)Ax\cr
\hfill&=\zero\nabla_{[\bar Y,\bar Z]}\bar T+{2\over\mu_{0}}\omega(\bar 
Y,\bar Z)\zero\nabla_{Ax}\bar T\cr
\hfill&\kern12pt+{1\over\mu_{0}}\omega(\bar T,A[\bar Y,\bar 
Z]+{2\over\mu_{0}}\omega(\bar Y,\bar Z)A^{2}x)x\cr
\hfill&\kern12pt+{1\over\mu_{0}}\omega([\bar Y,\bar 
Z]+{2\over\mu_{0}}\omega(\bar Y,\bar Z)Ax,\bar T)Ax.}$$
Hence:
$$\eqalign{\overline{R^{r}(Y,Z)T}&={1\over\mu_{0}}\omega(\bar T,A\bar 
Z)\bar Y-{1\over\mu_{0}}\omega(\bar T,A\bar Y)\bar 
Z+{1\over\mu_{0}}\omega(\bar Z,\bar T)A\bar Y\cr
\hfill&\kern12pt-{1\over\mu_{0}}\omega(\bar Y,\bar T)A\bar 
Z-{2\over\mu_{0}}\omega(\bar Y.\bar Z)A\bar T\cr
\hfill&\kern12pt+{1\over\mu_{0}^{2}}[\omega(\bar Z,\bar T)\omega(\bar 
Y,A^{2}x)-\omega(\bar Y,\bar T)\omega(\bar Z,A^{2}x)\cr
\hfill&\kern38pt-2\omega(\bar Y,\bar Z)\omega(\bar T,A^{2}x)]x\cr
\hfill&={1\over\mu_{0}}[\omega(\bar Z,\bar T)(A\bar 
Y-{1\over\mu_{0}}\omega(A\bar Y,Ax)x)\cr
\hfill&\kern28pt-\omega(\bar Y,\bar T)(A\bar 
Z-{1\over\mu_{0}}\omega(A\bar Z,Ax)x)\cr
\hfill&\kern28pt-2\omega(\bar Y,\bar Z)(A\bar 
T-{1\over\mu_{0}}\omega(A\bar T,Ax)x)\cr
\hfill&\kern28pt+\omega(\bar T,A\bar Z)\bar 
Y-\omega(\bar T,A\bar Y)\bar Z].\cr}$$
Notice that
$$A\bar Z-{1\over\mu_{0}}\omega(A\bar Z,Ax)x\an H_{x}$$
$$[A\bar Z-{1\over\mu_{0}}\omega(A\bar Z,Ax)x,Ax]=0$$
and thus we can define an operator $A_{y}$ on $T_{y}M$ by
$$\overline{A_{y}U}=A\bar U-{1\over\mu_{0}}\omega(A\bar 
U,Ax)x,\quad\hbox{where}\quad y=\pi(x).$$
Observe that
$$\eqalign{\Omega_{y}(A_{y}U,V)&=\omega(A\bar U-{1\over\mu_{0}}\omega(A\bar 
U,Ax)x,\bar V)\cr
\hfill&=\omega(A\bar U,\bar V)=-\omega(\bar U,A\bar 
V-{1\over\mu_{0}}\omega(A\bar V,Ax)x)\cr
\hfill&=-\Omega_{y}(U,A_{y}V),\cr}$$
i.e. $A_{y}$ is an element of the symplectic Lie algebra at $y\an M$. 
We can thus rewrite the curvature formula as:
$$\eqalign{R^{r}_{y}(Y,Z)T&={1\over\mu_{0}}[\Omega(Z,T)A_{y}Y-\Omega(Y,T)A_{y}Z\cr
\hfill&\kern28pt-2\Omega(Y,Z)A_{y}T+\Omega(T,A_{y}Z)Y-\Omega(T,A_{y}Y)Z].\cr}$$
%Hence
%$${1\over n+1}\Omega(X,AY)=-{2\over\mu_{0}}\Omega(X,\tilde AY)={1\over 
%n+1}Ric^{r}(X,Y),$$
which proves the theorem.\qed

\math{Theorem.}{\Asquare}{\sl The 
connection $\nabla^{r}$ is locally symmetric if and only if 
$A^{2}=\lambda\;{\rm Id}$.}

\vskip0.3cm

\undertext{\it Proof.} This is a  corollary of Theorem $\redcurv$. 
%In particular, one 
%sees that the covariant derivative of the Ricci tensor is given by:
%$$(\nabla^{r}_{W}Ric^{r})(U,V)={1\over 
%2n+1}[\Omega(W,U)u(V)+\Omega(W,V)u(U)]$$
%where
%$$u_{y=\pi(x)}(V)=-{2(n+1)(2n+1)\over\mu_{0}^{2}}
%\Omega_{y}(\pi_{\ast}(A^{2}x-{1\over\mu_{0}}\omega(A^{2}x,Ax)x),V).$$
Indeed, from the above one sees that
$$Ric^{r}(U,V)=\kappa\Omega(U,AV)$$
where $\kappa$ is a constant and $Ric^{r}$ is the Ricci curvature of 
the reduced connection. Hence:
$$\eqalign{(\nabla^{r}_{X}Ric^{r})(U,V)&=\kappa[X\Omega(U,AV)-
\Omega(\nabla^{r}_{X}U,AV)-\Omega(U,A\nabla^{r}_{X}V)]\cr
\hfill&=\kappa[\bar X\omega(\bar 
U,\overline{AV})-\omega(\overline{\nabla^{r}_{X}U},\overline{AV})-
\omega(\bar U,\overline{A\nabla^{r}_{X}V})]\cr
\hfill&=\kappa[\omega(\zero\nabla_{\bar X}\bar 
U,\overline{AV})+\omega(\bar U,\zero\nabla_{\bar 
X}\overline{AV})-\omega(\zero\nabla_{\bar X}\bar U,\overline{AV})\cr
\hfill&\kern20pt-\omega(\bar U,A(\zero\nabla_{\bar X}\bar 
V+{1\over\mu_{0}}\omega(\bar V,A\bar X)x+{1\over\mu_{0}}\omega(\bar 
X,\bar V)Ax))]\cr
\hfill&=\kappa[\omega(\bar U,\zero\nabla_{\bar X}(A\bar 
V+{1\over\mu_{0}}\omega(\bar V,A^{2}x)x))\cr
\hfill&\kern20pt-\omega(\bar 
U,A\zero\nabla_{\bar X}\bar V+{1\over\mu_{0}}\omega(\bar X,\bar 
V)A^{2}x)]\cr
\hfill&={\kappa\over\mu_{0}}[\omega(\bar V,A^{2}x)\omega(\bar U,\bar 
X)+\omega(\bar U,A^{2}x)\omega(\bar V,\bar X)].\cr}$$
Define the 1-form $u$ on $M$ by:
$$u_{y}(V)={\kappa\over\mu_{0}}\omega(\bar V,A^{2}x),\;y=\pi(x).$$
This has a meaning as
$$Ax[\omega(\bar V,A^{2}x)]=0.$$
The condition for local symmetry is
$$\omega(\bar V,A^{2}x)=0,\;\forall V$$
and it is identically satisfied if $A^{2}=\lambda\;{\rm Id}$. The 
converse is an immediate consequence of Theorem 2 of {\cahen}.\qed

\vskip0.3cm

\math{Corollary.}{\notsymmetric}{\sl Let $A\an\frak s\frak 
p(n+1,\R)$, $A^{2}\neq\lambda\;{\rm Id}$ and assume that the 
projection $\Sigma\rightarrow M$ is a smooth submersion. Then, $M$ 
admits a symplectic connection with curvature of Ricci type which is 
not locally symmetric.}

\vskip0.3cm

\math{Theorem.}{\symmspaces}{\sl The reduced spaces $M$ obtained 
above, when $A^{2}=\lambda\;{\rm Id}$, are globally symmetric.}

\undertext{\it Proof.} The argument (identical in all cases 
$\lambda>0,\lambda<0,\lambda=0$) is that the group $G$ of symplectic 
transformations of $\R^{2n+2}$ which commute with $A$, acts 
transitively on the quadric and this action projects onto an action 
of $G$ on $M$. Furthermore, the isotropy subgroup of a point of $M$ is 
the group of fixed points of an involutive automorphism of $G$. 
Finally the action of $G$ on $M$ is symplectic and affine.\qed

\vskip0.5cm

{\bf 5.} The extreme simplicity of this construction suggests to 
generalize it to the situation where one has a family of elements 
$A_{j}$ ($j\leq p$) of the symplectic Lie algebra with 
$[A_{i},A_{j}]=0,\;\forall i,j$. The submanifold
$$\omega(x,A_{j}x)=\mu_{j}\neq 0$$
can be reduced under the action of $\R^{p}$. In generic situation, 
the reduced space $M$ is indeed a differentiable manifold and a reduced 
connection may be defined in a similar way.

\vfill\eject
\vskip1cm
   \ifreferenceopen \Closeout\referencewrite \referenceopenfalse \fi
   \line{\bf\hskip0pt\hfil References\hfil}\vskip\headskip
   \vskip0.3cm
   \input referenc.txa
   %\footline={\hss\eightsc\folio\hss}\rm

\end